\documentclass[12pt,reqno]{article}
\usepackage{dcolumn}
\usepackage{amsthm,epsfig}
\usepackage{array,dcolumn,graphicx,amsmath}
\usepackage{booktabs,lscape,multirow,psfrag,rotating}
\usepackage{amssymb,amsmath,amsthm}
\usepackage{longtable}
\usepackage{booktabs,lscape,multirow,rotating}
\usepackage[normal]{caption2}
\usepackage[mathscr]{eucal}
\usepackage{eucal}
\usepackage{amssymb,pstricks}
\newcommand{\bSigma}{\boldmath$\Sigma$}
\newcommand{\mbSigma}{\mbox{\bSigma}}
\newcommand{\bTheta}{\boldmath$\Theta$}
\newcommand{\mbTheta}{\mbox{\bTheta}}

\newcommand{\bGamma}{\boldmath$\Gamma$}
\newcommand{\mbGamma}{\mbox{\bGamma}}

\newcolumntype{d}[1]{D{.}{.}{#1}}
\begin{document}
\def\s#1{\oalign{$#1$\crcr\hidewidth \normal$\sim$ \hidewidth}}
\thispagestyle{empty}
\newtheorem{theorem}{\indent THEOREM}
\newtheorem{prop}{\indent Proposition}
\newtheorem{lemma}{\indent Lemma}
\renewcommand{\proofname}{\hspace*{\parindent}{Proof.}}

\begin{center}
{\bf A note on MLE of covariance matrix}
\end{center}

\vspace{0.3cm}
\begin{center}
\parbox{13cm}{\sloppy \, \,
\noindent  {Ming-Tien Tsai}
 }
\end{center}

\begin{center}
\parbox{13cm}{\sloppy \, \,
\noindent  {Institute of Statistical Science, Academia Sinica, Taipei, Taiwan}

 }
\end{center}

\vspace{0.3cm}
\begin{center}
\parbox{13cm}{\sloppy \, \,
{\bf Summary:} For a multivariate normal set up, it is well known that the maximum likelihood estimator of
covariance matrix is neither admissible nor minimax under the Stein loss function. For the past six decades,
a bunch of researches have followed along this line for Stein's phenomenon in the literature. In this note,
the results are two folds: Firstly, with respect to Stein type loss function we use the full Iwasawa
decomposition to enhance the unpleasant phenomenon that the minimum risks of maximum likelihood estimators
for the different coordinate systems (Cholesky decomposition and full Iwasawa decomposition) are different.
Secondly, we introduce a new class of loss functions to show that the minimum risks of maximum likelihood
estimators for the different coordinate systems, the Cholesky decomposition and the full Iwasawa decomposition,
are of the same, and hence the Stein's paradox disappears.
}

\vspace{0.3cm}
\parbox{13cm}{\sloppy \, \,
{\it Keywords:} Geodesic distance; Iwasawa decomposition; minimax estimator}
\end{center}

\vspace{0.3cm}
\def \theequation{\arabic{equation}}
\setcounter{equation}{0}
\noindent {\bf 1. Introduction}
\vspace{0.3cm}

\indent  Let ${\bf X}_{1}, \cdots, {\bf X}_{n}$ be independent $p$-dimensional random vectors with a common
multivariate normal distribution $N_{p}(\bf 0, {\mbSigma})$. A basic problem considered in the literature
is the estimation of the $p \times p$ covariance matrix ${\mbSigma}$ which is unknown and assumed to be
nonsingular. It is also assumed that $n \geq p$, as such the sufficient statistic
\begin{align}
{\bf A}=\sum_{i=1}^{n} {\bf X}_{i}{\bf X}^{'}_{i}
\end{align}
is positive definite with probability one. In the literature, the estimators $\phi({\bf A})$ of ${\mbSigma}$
are the functions of ${\bf A}$. The sample space ${\mathcal S}$, the parameter space ${\varTheta}$ and the
action space ${\mathcal A}$ are taken to be the set of $p \times p$ positive definite matrices. Note that
${\bf A}$ has a Wishart distribution $W({\mbSigma}, n)$ and the maximum likelihood estimator of ${\mbSigma}$
\begin{align}
\Hat{\mbSigma}=n^{-1}{\bf A},
\end{align}
which is unbiased. The general linear group $Gl(p)$ acts on the spaces ${\mathcal S}$, ${\varTheta}$ and
${\mathcal A}$. Generally, we consider the invariance loss function $L$, namely, $L$ satisfies the condition
that $L(g\phi({\bf A})g^{'}, g{\mbSigma}g^{'})=L(\phi({\bf A}), {\mbSigma})$ for all $g \in Gl(p)$.  One of
the most interesting examples was introduced by Stein (see Jame and Stein, 1961),
\begin{align}
L(\phi({\bf A}), {\mbSigma})=\mbox{tr}{\mbSigma}^{-1}\phi({\bf A})
      -\mbox{log}\mbox{det}{\mbSigma}^{-1}\phi({\bf A})-p,
\end{align}
where $\mbox{tr}$ and $\mbox{det}$ denote the trace and the determinant of a matrix, respectively.
Because $Gl(p)$ acts transitively on the space ${\varTheta}$, so the best equivalent estimator exists.
The minimum risk for the estimator $\hat{\mbSigma}$ is
\begin{align}
R(\Hat {\mbSigma}, {\mbSigma})=\sum_{i=1}^{p}\{ \mbox{log} n
   - {\mathcal E}[\mbox{log}{\chi}^{2}_{n-i+1}]\},
\end{align}
where ${\mathcal E}[X]$ denotes the expectation of random variable of $X$. It can be easily seen that the
maximum likelihood estimator is the best equivalent estimator.




\indent Since the general linear group is not a solvable group, hence relax the condition a little bit by
considering the group of $p \times p$ lower triangular matrices with positive diagonal elements $G^{+}_{T}$,
the loss function is also invariant under $G^{+}_{T}$. Using the Cholesky decomposition, we may write
${\bf A}={\bf T}{\bf T}^{'}$, where ${\bf T} \in G^{+}_{T}$. Since $G^{+}_{T}$ acts
transitively on the space ${\varTheta}$, the best equivalent estimator was established by Stein
(see James and Stein, 1961) in the following
\begin{align}
\Hat{\mbSigma}_{S}={\bf T}{\bf D}^{-1}_{S}{\bf T}^{'},
\end{align}
where $D_{S}$ is a positive diagonal matrix with elements $d_{Sii}=n+p-2i+1, ~i=1, \cdots, p$. The minimum risk
for the estimator $\hat{\mbSigma}_{S}$ is
\begin{align}
R(\Hat {\mbSigma}_{S}, {\mbSigma})=\sum_{i=1}^{p}\{\mbox{log} (n+p-2i+1)
           -{\mathcal E}[\mbox{log}{\chi}^{2}_{n-i+1}]\}.
\end{align}
Since the group $G^{+}_{T}$ is solvable, it follows from results in Kiefer (1957) that the estimator
$\hat{\mbSigma}_{JS}$ is minimax.

\indent In the literature, there are many developements along this approach and its ramplifications, we may refer
to the book of Anderson (2003) or the book of Muirhead (1982), and the references cited there, hence we omit the
details. With respect to Stein loss function, we use the full Iwasawa decomposition (Terras, 1988) to enchance the
Stein's phenomenon.

\vspace{0.3cm}
\noindent {\bf 2. The full Iwasawa decomposition}
\vspace{0.3cm}

\indent The Cholesky decomposition can be viewed as a partial Iwasawa decomposition. We would like to relax
the conditions more by considering the full Iwasawa decomposition in this section. Some more notations are
needed. Partition ${\mbSigma}_{(k)}$ and ${\bf A}_{(k)}$ as
\begin{eqnarray}
  {\mbSigma}_{(k)} = \left[
  \begin{array}{cc}
    {\sigma}_{(k)11} & {\mbSigma}_{(k)12}\\
    {\mbSigma}_{(k)21} & {\mbSigma}_{(k)22}
  \end{array}
  \right]
  ~~\mbox{and}~~
   {\bf A}_{(k)} = \left[
  \begin{array}{cc}
    {a}_{(k)11} & {\bf A}_{(k)12}\\
    {\bf A}_{(k)21} & {\bf A}_{(k)22}
  \end{array}
  \right],
\end{eqnarray}
for all $k=1, \cdots, p$ with ${\mbSigma}_{(1)}={\mbSigma}$ and ${\bf A}_{(1)}={\bf A}$, also define
\begin{align}
{\mbSigma}_{(k+1)}=\mbSigma_{(k)22}- \mbSigma_{(k)21} \mbSigma_{(k)12}/\sigma_{(k)11}
\end{align}
and
\begin{align}
{\bf A}_{(k+1)}={\bf A}_{(k)22}- {\bf A}_{(k)21}{\bf A}_{(k)12}/a_{(k)11}.
\end{align}
Let
\begin{align}
{\bf g}_{(k)} = \left[
\begin{array}{cc}
 1  & {\bf 0}    \\
-\mbSigma_{(k)21}\sigma_{(k)11}^{-1} & {\bf I}
\end{array}
\right] ~~\mbox{and}~~
{\bf h}_{(k)} = \left[
\begin{array}{cc}
  1  & {\bf 0}   \\
-{\bf A}_{(k)21}a_{(k)11}^{-1} & {\bf I}
\end{array}
\right], ~k=1, \cdots, p.
\end{align}
Then we have
\begin{align}
\widetilde{\mbSigma}_{(k)} &= {\bf g}_{(k)} {\mbSigma}_{(k)} {\bf g}^{'}_{(k)}   \\  \nonumber
&= \left[
\begin{array}{cc}
\sigma_{(k)11}& {\bf 0}\\
{\bf 0}& \mbSigma_{(k)22:1}
\end{array}
\right],
\end{align}
and
\begin{align}
\widetilde{\bf A}_{(k)} &= {\bf h}_{k} {\bf A}_{(k)} {\bf h'}_{(k)}= \left[
\begin{array}{cc}
 a_{(k)11}& {\bf 0}\\
{\bf 0}& {\bf A}_{(k)22:1}
\end{array}
\right], ~~k=1, \cdots, p.
\end{align}
Let
\begin{align}
{\mbSigma}^{*}=\mbox {Diag}({\sigma}_{(1)11}, \cdots, {\sigma}_{(p)11}) ~ \mbox{and}~
{\bf A}^{*}=\mbox {Diag}(a_{(1)11}, \cdots, a_{(p)11}).
\end{align}
By using the full Iwasawa decomposition, we can eventually transform ${\mbSigma}$ and ${\bf A}$ into the
diagonal matrices ${\mbSigma}^{*}$ and ${\bf A}^{*}$, respectively. Thus we establish the one-to-one
correspondences: ${\mbSigma} \leftrightarrow {\mbSigma}^{*}$ and ${\bf A} \leftrightarrow {\bf A}^{*}$. By the
properties of Wishart distribution (see Theorem 4.3.4, Theorem 7.3.4 and Theorem 7.3.6 of Anderson, 2003),
it is easy to note that $a_{(i)11}/{\sigma}_{(i)11}, i=1, \cdots, p$ are independent $\chi^{2}$ random variables
with $n-i+1$ degrees of freeedom respectively. Consider the loss function similar to the equation (3)
\begin{align}
  L(\phi ({\bf A}^{*}), {\mbSigma}^{*})=\mbox{tr}{\mbSigma}^{*-1}{\bf D}{\bf A}^{*}
     -\mbox{log}\mbox{det}{\mbSigma}^{*-1}{\bf D}{\bf A}^{*}-p,
\end{align}
where ${\bf D}=\mbox {Diag}(d_{11}, \cdots, d_{pp})$ is a positive diagonal matrix, not depending on
${\bf A}^{*}$.

\vspace{0.3cm}
\indent {\bf Theorem 1.} {\it With respect to the likelihood loss function (Stein type loss function), the best
estimator invariant with respect to one-to-one transformations ${\mbSigma} \to {\mbSigma}^{*}, {\bf A}
\to {\bf A}^{*}$, is $\Hat{\mbSigma}^{*}_{I}={\bf D}^{-1}_{0}{\bf A}^{*}$. The minimum risk is
${\mathcal E} L(\phi ({\bf A}^{*}), {\bf I})=\sum_{i=1}^{p}\{\mbox{log} (n-i+1)
-{\mathcal E}[\mbox{log}{\chi}^{2}_{n-i+1}]$, and the estimator is minimax.}
\vspace{0.3cm}

\indent {\bf Proof}. Take ${\mbSigma}={\bf I}$, and then note that
\begin{align}
{\mathcal E} L(\phi ({\bf A}^{*}), {\bf I})&={\mathcal E}[\mbox{tr}{\bf D}{\bf A}^{*}
     -\mbox{log}\mbox{det}{\bf D}{\bf A}^{*}-p]   \\ \nonumber
     &=\sum_{i=1}^{p}(n-i+1)d_{ii}-\sum_{i=1}^{p}\mbox{log}d_{ii}
     -\sum_{i=1}^{p}{\mathcal E}[\mbox{log}{\chi}^{2}_{n-i+1}]-p.
\end{align}
The minimum of (15) occurs at $d_{ii}=1/(n-i+1), i=1, \cdots, p$. Since ${\bf A}^{*}$ also acts transitively on
the space ${\varTheta}$, so the best equivalent estimator exists, which is of the form
\begin{align}
  \Hat{\mbSigma}^{*}_{I}={\bf D}^{-1}_{0}{\bf A}^{*},
\end{align}
where ${\bf D}_{0}$ is a diagonal matrix with elements $d_{0ii}=n-i+1, i=1, \cdots, p$.
Thus the minimum risk for the estimator $\hat{\mbSigma}^{*}_{I}$ is
\begin{align}
  R(\Hat {\mbSigma}^{*}_{I}, {\mbSigma}^{*})=\sum_{i=1}^{p}\{\mbox{log} (n-i+1)
           -{\mathcal E}[\mbox{log}{\chi}^{2}_{n-i+1}]\}.
\end{align}

\indent Since the group of positive diagonal matrices is a subset of $G^{+}_{T}$, which is solvable, thus the
group of positive diagonal matrices is also solvable. And hence, by the results of Kiefer (1957) the estimator
$\hat{\mbSigma}^{*}_{I}$ is minimax.

\indent Compare equations (4), (6) and (17), then it is easily to see that
\begin{align}
 R(\Hat {\mbSigma}^{*}_{I}, {\mbSigma}^{*}) \leq R(\Hat {\mbSigma}_{S}, {\mbSigma})
        \leq R(\Hat {\mbSigma}, {\mbSigma}),
\end{align}
for $p \geq 2$. 
The equality in (18) holds  when $\Hat {\mbSigma}^{*}_{I}=\Hat {\mbSigma}_{S}=\Hat {\mbSigma}$:
namely, (i) the components are independent or (ii) as the sample size $n \to \infty$. With respect to the Stein
loss function, the minimum risk functions are different based on the full Iwasawa decomposition and based
on the Cholesky decomposition.

\indent We may note that each $a_{(i)11}/(n-i+1)$ is the maximum likelihood estimator of ${\sigma}_{(i)11}$,
and is unbiased, for all $i=1, \cdots, p$. For each component, $a_{(i)11}/(n-i+1)$ is admissible for
${\sigma}_{(i)11}, i=1, \cdots, p$. Note that saying $n^{-1}{\bf A}$ is the maximum likelihood estimator of
${\mbSigma}$ is the same as to say that ${\bf D}^{-1}_{0}{\bf A}^{*}$ is the maximum likelihood estimator of
${\mbSigma}^{*}$. Thus, the results of the equation (18) lead to a paradox that the property of maximum likelihood
estimators for the different coordinate systems is not consistent with respect to  the Stein type loss function.
This motives us to further study whether a suitable loss function exists so that the maximum likelihood estimators
can be invariant under reparameterizations.


\vspace{0.3cm}
\noindent {\bf 3. The geodesic distance}
\vspace{0.3cm}

\indent Since the space of positive definite symmetric matrices is a non-Euclidean space, it is more natural
to use a metric on a Riemannian metric space. The Riemannian metric can be defined with the help of the
fundamental form $ds^{2}= \mbox {tr} ({\bf W}^{-1}d{\bf W})^{2},$ where $d{\bf W}$ denotes the matrix of
differentials. This is invariant under the transformation ${\bf W} \to {\bf V}{\bf W}$, where ${\bf V}$ is any
fixed elements of $Gl(p)$. Let $\mathcal P_{p}$ be the set of square symmetric positive definite matrices,
this set is a representation space of the group $Gl(p)$. An element ${\bf V} \in Gl(p)$ operates on
${\mathcal P}_{p}$ according to ${\bf M} \to {\bf V}{\bf M}{\bf V}^{'}$. On $Gl(p)$, any maximal compact
subgroup $\triangle$ of $Gl(p)$ can be represented in the form
$\triangle=\{{\bf V}| {\bf V}^{'}{\bf H}{\bf V}={\bf H}, {\bf V} \in Gl(p)\}$. The set of cosets
$Gl(p)/\triangle=\{\rho\triangle| ~\rho \in Gl(p)\}$ can be considered as a representation space of
$Gl(p)$. Since $\triangle={\bf V}^{-1}{\mathcal O}(p){\bf V}, {\bf V} \in Gl(p)$, thus the all maximal compact
subgroup are conjugate. Conjugate subgroups yield equivalent representation spaces, hence it is sufficiently
enough to only consider the orthogonal group ${\mathcal O}(p)$ as the maximal compact group of $Gl(p)$.
The map ${\bf V}{\mathcal O}(p) \to {\bf M}={\bf W}{\bf W}^{'}$ establishes an equivalence between the
representation spaces $Gl(p)/{\mathcal O}(p)$ and $\mathcal P_{p}$ of $Gl(p)$. And hence the $ds^{2}$ defines
an invariant metric on ${\mathcal P}_{p}$. The tangent space to ${\mathcal O}(p)$ is $so(p)$, the vector
space of skew-symmetric $p \times p$ matrices. Thus $\mbox{dim}{\mathcal O}(p)=\mbox{dim}so(p)=p(p-1)/2$.
And hence $\mbox{dim}Gl(p)/{\mathcal O}(p)=\mbox{dim}Gl(p)-\mbox{dim}{\mathcal O}(p)=p(p+1)/2-p(p-1)/2=p$.

\vspace{0.3cm}
\indent {\bf Proposition 1.}  {\it A geodesic segment $T(t)$ through ${\bf I}$ and ${\bf Y}$ in ${\mathcal P}_{p}$
has the form: $T(t)=\mbox{exp}(t{\bf U}^{'}{\bf B}{\bf U}), 0 \leq t \leq 1$, where ${\bf Y}$ has the spectral
decomposition: ${\bf Y}={\bf U}^{'}\mbox{exp}({\bf B}){\bf U}=\mbox{exp}({\bf U}^{'}{\bf B}{\bf U})$,
for ${\bf U} \in {\mathcal O}(p)$ and ${\bf B}=\mbox{Diag}(b_{1}, \cdots, b_{p}), b_{i} \in R, i=1, \cdots, p.$
The length of the geodesic segment is $(\sum_{i=1}^{p}b^{2}_{i})^{1/2}.$}
\vspace{0.3cm}

\indent The proof of Proposition 1 can be found in the book of Terras (1988). For any given two points
$\mbSigma_{1}$ and $\mbSigma$ of $\mathcal P_{p}$, the geodesic distance is defined to be
$\sum_{i=1}^{p}\mbox{log}^{2}\lambda_{i}$, where $\lambda_{i}, i=1, \cdots, p,$ are the zeros of charactistic
polynomial $\mbox{det}(\lambda{\mbSigma}-{\mbSigma}_{1})$.  A loss function $L$ is invariant iff $L$ can
be written as a function of the eigenvalues of $\phi({\bf A})$. Hence from geometric point of view, with
the help of Proposition 1 we may naturally consider the compatible, coordinate-invariant loss function
$L_{G}(\phi({\bf A}), {\mbSigma})=\sum_{i=1}^{p}\mbox{log}^{2}\lambda_{i}$, where $\lambda_{i}, i=1, \cdots, p$
denote the eigenvalues of ${\mbSigma}^{-1}\phi({\bf A})$. And the risk function is denoted by
$R_{G}(\Hat {\mbSigma}_{G}, {{\mbSigma}})(={\mathcal E}[L_{G}(\phi({\bf A}), {\mbSigma})])$, where
$\Hat{\mbSigma}_{G}$ be the corresponding estimator based on the geodesic distance loss function on
$\mathcal P_{p}$. Without loss the generality we may take ${\mbSigma}={\bf I}$ as Stein did.
Write ${\bf d}=(d_1, \cdots, d_p)^{'}$. The following studies are comparable to the results in Section 7.8.
of Anderson (2003). Let ${\mathcal Var}[X]$ denote the variance of random variable $X$.

\vspace{0.3cm}
\indent {\bf Theorem 2.} {\it With respect to the geodesic distance  loss function on $\mathcal P_{p}$,
$L_{G}(\phi({\bf A}), {\bf I})=\sum_{i=1}^{p}\mbox{log}^{2}(d_{i}\lambda_{i})$, where $\lambda_{i}$ denotes
the $i$-th largest eigenvalue of $\phi({\bf A})$ and $d_{i}$ is a positive constant, $\forall i=1, \cdots, p$,
let ${\mathcal C}=\{{\bf d}|{\mathcal E}[d_{i}\lambda_{i}]\leq e, \forall i=1, \cdots, p\}$.
Then on the set ${\mathcal C}$ the minimum of risk function $R_{G}(\Hat {\mbSigma}_{G}, {\bf I})$ occurs at
$d_{i}=\mbox{exp}\{-{\mathcal E}[\mbox{log}\lambda_{i}]\}, \forall i=1, \cdots, p$, and its minimum risk is
$\sum_{i=1}^{p}{\mathcal Var}[\mbox{log}\lambda_{i}]$.}
\vspace{0.3cm}

\indent {\bf Proof}. Note that $\partial {\mathcal E}[\sum_{i=1}^{p}\mbox{log}^{2}(d_i\lambda_{i})]
/\partial {d_{i}}=0$ implies that ${\mathcal E}[\mbox{log}(d_i\lambda_{i})/d_{i}]=0, \forall i=1, \dots, p$,
which is the same as the conditions that $\mbox{log}d_i+{\mathcal E}[\mbox{log}\lambda_{i}]=0,
\forall i=1, \cdots, p$. By Jensen inequality, we have that
$0=\mbox{log}d_i+{\mathcal E}[\mbox{log}\lambda_{i}], \forall i=1, \cdots, p$. Let ${\bf d}_{0}$ be the
point so that $\mbox{log}d_{i}=-{\mathcal E}[\mbox{log}\lambda_{i}], \forall i=1, \cdots, p$.
Then ${\bf d}_{0}$ is the critical point for the risk function $R_{G}(\Hat {\mbSigma}_{G}, {\bf I})$.

\indent Moreover, we may note that $\partial^{2} {\mathcal E}[\sum_{i=1}^{p}\mbox{log}^{2}(d_i\lambda_{i})]
/\partial^{2} {d_{i}}=2{\mathcal E}[(1-\mbox{log}d_i\lambda_{i})/d^{2}_{i}], \forall i=1, \cdots, p$. By the
Jensen inequality, we may note that ${\mathcal E}[(1-\mbox{log}d_i\lambda_{i})] \geq
(1-\mbox{log}{\mathcal E}[d_i\lambda_{i}]) \geq 0, \forall i=1, \cdots, p$. Then, on the set
${\mathcal C}$, we have $\partial^{2}{\mathcal E}[\sum_{i=1}^{p}\mbox{log}^{2}(d_i\lambda_{i})]
/\partial^{2} {d_{i}} \geq 0, \forall i=1, \cdots, p$, and
$\partial^{2} {\mathcal E}[\sum_{i=1}^{p}\mbox{log}^{2}(d_i\lambda_{i})]
/\partial d_{i}\partial d_{j}=0, \forall i\ne j$. Thus the risk function $R_{G}(\Hat {\mbSigma}_{G}, {\bf I})$
is convex and has an unique minimum on the set ${\mathcal C}$. Since the set ${\mathcal C}$ is a connected set,
and ${\bf d}_{0} \in {\mathcal C}$, and hence, the theorem follows.
\vspace{0.3cm}

\indent Theorem 2 provides us a new class of loss functions for statistical inference. Although Theorem 2
looks simple mathematically, it provides us an intrinsicaly new approach to make statistical inference.
With respect to the geodesic distance loss function, we may see that MLE is invariant under different
parameterizations, and more importantly, the Stein's paradox disappears. Those results are tremendously
different from what have existed in the literature with respect to Stein type loss functions.

\indent To obtain the best equivalent estimator with respect to geodesic distance loss function, the quantities
$\mbox{exp}\{-{\mathcal E}[\mbox{log}\lambda_{i}]\}, \forall i=1, \cdots, p$, which can be viewed as the
geometric means of the distributions, are needed to be evaluated. Some often seen cases are illustrated in
the followings:

\vspace {0.3cm}
\indent {\bf I. The full Iwasawa decomposition form}. Based on the one-to-one transformation:
${\mbSigma} \leftrightarrow {\mbSigma}^{*}$ and ${\bf A} \leftrightarrow {\bf A}^{*}$, as such
using the geodesic distance on $\mathcal P_{p}$, namely,  to use the loss function
$L_{G}({\bf D}{\bf A}^{*}, {\bf I})=\sum_{i=1}^{p}\mbox{log}^{2}(d_{ii}\lambda^{*}_{i})$, where
$\lambda^{*}_{i}, i=1, \cdots, p$ are the eigenvalues of ${\bf A}^{*}$. Thus,
${\mathcal E}[L_{G}({\bf D}{\bf A}^{*}, {\bf I})]
={\mathcal E}[\sum_{i=1}^{p}\mbox{log}^{2}(d_{ii}\lambda^{*}_{i})]$. Thus, by Theorem 2 the minimum of it
occurs at $d_{ii}=\mbox{exp}\{-{\mathcal E}[\mbox{log}{a_{(i)11}}]\}
=\mbox{exp}\{-{\mathcal E}[\mbox{log}{\chi}^{2}_{n-i+1}]\}, \forall i=1, \cdots, p$.
Thus the geometric estimator $\Hat{\mbSigma}^{*}_{GI}$ of ${\mbSigma}^{*}$ is that
$\Hat{\mbSigma}^{*}_{GI}={\bf D}^{*}{\bf A}^{*}$, where ${\bf D}^{*}$ is a diagonal matrix with elements
$d^{*}_{ii}=\mbox{exp}\{-{\mathcal E}[\mbox{log}{\chi}^{2}_{n-i+1}]\}, i=1, \cdots, p$ and ${\bf A}^{*}$ in (13) .
Since ${\bf A}^{*}$ is a diagonal matrix, the group of positive diagonal matrices is solvable, and hence by
the results of Kiefer (1957) this geometric estimator $\Hat{\mbSigma}^{*}_{GI}$ is minimax. And its  minimum
risk is $R_{G}(\Hat {\mbSigma}^{*}_{GI}, {\bf I})=\sum_{i=1}^{p}{\mathcal Var}[\mbox{log}{\chi}^{2}_{n-i+1}]$.
Also note that $R_{G}(\Hat {\mbSigma}^{*}_{I}, {\bf I})=R_{G}(\Hat {\mbSigma}^{*}_{GI}, {\bf I})
+\sum_{i=1}^{p}\{\mbox{log}(n-i+1)-{\mathcal E}[\mbox{log}{\chi}^{2}_{n-i+1}]\}^{2}, i.e.,
R_{G}(\Hat {\mbSigma}^{*}_{GI}, {\bf I}) < R_{G}(\Hat {\mbSigma}^{*}_{I}, {\bf I})$, where
$\Hat {\mbSigma}^{*}_{I}$ is defined in (16). Thus we may conclude that the Stein estimator
$\Hat {\mbSigma}^{*}_{I}$ is inadmissible with respect to geodesic distance loss function.

\vspace{0.3cm}
\indent {\bf II. The Cholesky decomposition form, can be viewed as a partial Iwasawa decomposition}.
Riemannian geometry yields an invariant volume element $dv$ on $\mathcal P_{p}$, it is of the form
$dv=(\mbox{det}{\bf A})^{-(p+1)/2}(d{\bf A})$, where $(d{\bf A})=\prod_{1 \leq j \leq i \leq p}d{a}_{ij}$ with
$d{a}_{ij}$ being the Lebesgue measure on $R^{p(p+1)/2}$. This $dv$ is the invariant $d$-form $(d=p(p+1)/2)$
on $\mathcal P_{p}$. Normalization is not necessary because this invariant $d$-form is not a probability measure.
Note that this $d$-form is still invariant under $G^{+}_{T}$. With the Cholesky decomposition
${\bf A}={\bf T}{\bf T}^{'}$ and differentiate at ${\bf T}={\bf I}: d{\bf A}=d{\bf T}+d{\bf T}^{'}$.
Thus $da_{ii}=2 dt_{ii}$, and for $i > j, da_{ij}=dt_{ij}$. Then the $d$-form becomes to
$dv=2^{p}\prod_{i=1}^{p}t^{-i}_{ii}(dT)$, where $(dT)$ denotes the wedge product. Let
$\mbox{etr}({\bf A})$ denote $\mbox{exp}(\mbox{tr}{\bf A})$ and write ${\bf T}=((t_{ij}))$.
Similarly, we do the Cholesky decomposition for the scale parameter as ${\mbSigma}={\mbTheta}{\mbTheta}^{'}$,
where ${\mbTheta}=((\theta_{ij}))$ is a lower triangular matrix. Take ${\mbSigma}={\bf I}$, the Wishart density
of ${\bf A}$ then becomes to the following density function
\begin{align}
\frac{2^{p-pn/2}}{\Gamma_{p}(n/2)}\mbox{etr}(-\frac{1}{2}{\bf T}{\bf T}^{'})\prod_{i=1}^{p}t^{n-i}_{ii},
\end{align}
where $\Gamma_{p}(a)=\pi^{p(p-1)/4}\prod_{i=1}^{p}\Gamma(a-(i-1)/2)$. This is essentially called the Bartlett
decomposition in the literature. Note that ${\bf T}$ is a lower triangular matrix which having the eigenvalues
$t_{ii}, i=1, \cdots, p$. With respect to  the geodesic distance loss function, the density function (19)
can be further reduced to the product of the $\chi^{2}$ density with $n-i+1$ degrees of freedom, $i=1, \cdots, p$.
Thus, the risk function of geodesic distance loss function is
${\mathcal E}[L_{G}(\phi({\bf A}), {\bf I})]={\mathcal E}[\sum_{i=1}^{p}\mbox{log}^{2}(d_{ii}t^{2}_{ii})]$.
By Theorem 2, the minimum of it occurs at $d_{ii}=\mbox{exp}\{-{\mathcal E}[\mbox{log}{t^{2}_{ii}}]\}
=\mbox{exp}\{-{\mathcal E}[\mbox{log}{\chi}^{2}_{n-i+1}]\}, \forall i=1, \cdots, p$. Let
$d^{0}_{ii}=\mbox{exp}\{-{\mathcal E}[\mbox{log}{t^{2}_{ii}}]\}, \forall i=1, \cdots, p$, and write
${\bf T}_{0}=((t_{0ij}))$, where $t_{0ij}=t_{ij}$ if $i\ne j$ and $t_{0ii}=(d^{0}_{ii}t^{2}_{ii})^{1/2},
i=1, \cdots, p, j=1, \cdots, i$. We may note that $t^{2}_{0ii}$ is the unbiased MLE of
${\theta}^{2}_{ii}, \forall i=1, \cdots, p$. Thus, the best equivalent estimator is of the
form $\Hat {\mbSigma}_{GC}={\bf T}_{0}{\bf T}^{'}_{0}$ and its minimum risk is
$R_{G}(\Hat {\mbSigma}_{GC}, {\bf I})=\sum_{i=1}^{p}{\mathcal Var}[\mbox{log}{\chi}^{2}_{n-i+1}]$.
Since the group $G^{+}_{T}$ is solvable, it follows from results in Kiefer (1957) that the estimator
$\Hat{\mbSigma}_{GC}$ is minimax with respect to geodesic distance loss function.

\indent This minimum risk is equivalent to that in Example 1, which indicates that with respect to the geodesic
distance loss function on the space of positive definite symmetric matrices the minimum risks of maximum
likelihood estimators with the different coordinate systems, the Cholesky decomposition and the full Iwasawa
decomposition, are of the same. These results are quite different from what have existed in the literature
based on the Stein type loss function, see Section 2 for the details. Moreover,  note that
$R_{G}(\Hat {\mbSigma}_{S}, {\bf I})=\sum_{i=1}^{p}{\mathcal E}[\mbox{log}^{2}({\chi}^{2}_{n-i+1}/n+p-2i+1)]$,
where $\Hat {\mbSigma}_{S}$ is defined in (5). It is easy to see that
$R_{G}(\Hat {\mbSigma}_{S}, {\bf I})=R_{G}(\Hat {\mbSigma}_{GC}, {\bf I})
+\sum_{i=1}^{p}\{\mbox{log}(n+p-2i+1)-{\mathcal E}[\mbox{log}{\chi}^{2}_{n-i+1}]\}^{2}, i.e.,
R_{G}(\Hat {\mbSigma}_{GC}, {\bf I}) < R_{G}(\Hat {\mbSigma}_{S}, {\bf I})$. Thus we may conclude
that the Stein estimator $\Hat {\mbSigma}_{S}$ is inadmissible with respect to geodesic distance loss function.

\vspace{0.3cm}
\indent {\bf III. The orthogonal decomposition form}. Stein (1956) considered the rotation-equivariant
estimator of $\mbSigma$. The class of rotation-equivariant estimators of covariance matrix is constituted
of all the estimators that have the same eigenvectors as the sample covariance matrix. Let $\lambda_i$ denotes
the $i$-th largest eigenvalue of ${\mbSigma}^{-1}{\bf A}$. Also write ${\bf A}={\bf U}{\bf L}{\bf U}^{'}$,
where ${\bf L}$ is a diagonal matrix with eigenvalues $l_{i}, i=1, \cdots, p$ and ${\bf U}$ being the
corresponding orthogonal matrix. Similarly, write ${\mbSigma}={\bf H}{\mbGamma}{\bf H}^{'}$,
where ${\mbGamma}$ is a diagonal matrix with eigenvalues $\gamma_{i}, i=1, \cdots, p$ and ${\bf H}$ being
the corresponding orthogonal matrix. Take ${\mbSigma}={\bf I}$, thus for the class of rotation-equivariant
estimators the minimum risk function based on the geodesic distance loss function is
${\mathcal E}[L_{G}(\phi({\bf A}), {\bf I})]={\mathcal E}[\sum_{i=1}^{p}\mbox{log}^{2}(d_{ii}l_{i})]$, and its
minimium occurs at $d_{ii}=\mbox{exp}\{-{\mathcal E}[\mbox{log}{l_{i}}]\}, \forall i=1, \cdots, p$.
Thus, the best rotation-equivariant estimator is of the form
$\Hat {\mbSigma}_{GO}={\bf U}{\bf D}^{*}{\bf L}{\bf U}^{'}$, where ${\bf D}^{*}$ is a diagonal matrix with
elements $d^{*}_{ii}=\mbox{exp}\{-{\mathcal E}[\mbox{log}{l_{i}}]\}, i=1, \cdots, p$. On the other hand,
the risk function based on the Stein loss function is ${\mathcal E} L(\phi ({\bf A}), {\bf I})
={\mathcal E}[\mbox{tr}{\bf D}{\bf L}-\mbox{log}\mbox{det}{\bf D}{\bf L}-p]=\sum_{i=1}^{p}d_{ii}
{\mathcal E}[l_{i}]-\sum_{i=1}^{p}\mbox{log}d_{ii}-\sum_{i=1}^{p}{\mathcal E}[\mbox{log}l_{i}]-p.$
The minimum of it occurs at $d^{-1}_{ii}={\mathcal E}[l_{i}], \forall i=1, \cdots, p$. Thus, with respect to
Stein loss function the best rotation-equivariant estimator for ${\mbSigma}$ is of the form
$\Hat {\mbSigma}_{O}={\bf U}{\bf D}^{*}_{0}{\bf L}{\bf U}^{'}$, where ${\bf D}^{*}_{0}$ is a diagonal matrix
with elements $d^{*}_{0ii}=1/{\mathcal E}[l_{i}], i=1, \cdots, p$. Similarly, we may also note that
$R_{G}(\Hat {\mbSigma}_{O}, {\bf I})=R_{G}(\Hat {\mbSigma}_{GO}, {\bf I})
+\sum_{i=1}^{p}\{\mbox{log}{\mathcal E}[l_{i}]-{\mathcal E}[\mbox{log}{l_{i}}]\}^{2}, i.e.,
R_{G}(\Hat {\mbSigma}_{GO}, {\bf I}) < R_{G}(\Hat {\mbSigma}_{O}, {\bf I})$. Thus, under the orthogonal
decomposition the Stein type estimator $\Hat {\mbSigma}_{O}$ is inadmissible with respect to geodesic distance
loss function.

\indent Via the result of Askey (1980), the joint density function of eigenvalues is of the form
\begin{align}
2^{-np/2}\prod_{i=1}^{p}\frac{\Gamma(3/2)}{\Gamma((1+i/2)\Gamma((n-p+i)/2)}
l^{\frac{n-p+i}{2}-1}_{i}\mbox{e}^{-{l_{i}}/2}\prod_{i < j}|l_i-l_j|.
\end{align}
It is easy to see that ${\mathcal E}[\prod_{i=1}^{p}\{l_i/(n-p+i)\}]=1$. However, we have difficulty to find
out the explicit form for the marginal density function of sample eigenvalue $l_i, i=1, \cdots, p$, and
open this type Selberg integral as a conjecture.

\vspace{0.3cm}
\noindent {\bf 4. General remarks}
\vspace{0.3cm}

\indent Entropy (expectation of likelihood loss function) not only plays an important role in information theory,
but also is a core in statistical theory. For the past six decades, quadratic loss function and likelihood loss
function are oftenly adopted to study the Stein's phenomenon for the covariance matrix estimation, for the details
see the Section 7.8. of Anderson (2003) or the Section 4.3. of Muirhead (1982). Via the full Iwasawa
decomposition, Theorem 1 tells us that the likelihood (Stein type) loss function is not invariant to arbitrary
reparameterizations of $\mathcal P_{p}$. To overcome the drawbacks, Riemannian metric is a natural way to be
adopted as a loss function for a non-Eculidean space $\mathcal P_{p}$. Due to the diffeomorphism invarinace of
risk function based on the geodesic distance, we may anticipate that the minimum risks of the MLEs may not only
be invariant to reparameterizations but also the Stein paradox disappear with respect to geodesic distance loss
function. Examples 1 and 2 tell us that the minimum risks of the MLEs of covariance matrices under the different
coordinate systems, the Cholesky decomposition and the full Iwasawa decomposition, are of the same with
respect to the geodesic distance loss function. Examples 1 and 2 also indicate that the MLE of covariance matrix
is minimax with respect to the geodesic distance loss function on $\mathcal P_{p}$. This note will inevitably have
an impact on statistical inference because the statisticans have to reconsider the adoptation of the MLE of
covariance matrix, which, however, had been constantly warned not to use for a long time since Stein's
phenonmenon occurred.

\indent If the quantity $d_{i}=\mbox{exp}\{-{\mathcal E}[\mbox{log}\lambda_{i}]\}$ in Theorem 2 can be
approximated by the quantity $d_{i}=1/{\mathcal E}[\lambda_{i}]$, $\forall i=1, \cdots, p$, then we have that
$\Hat {\mbSigma}^{*}_{GI}=\Hat {\mbSigma}^{*}_{I}$ and $\Hat {\mbSigma}_{GO}=\Hat {\mbSigma}_{O}$ approximatly.
Moreover, $R_{G}(\Hat {\mbSigma}^{*}_{GI}, {\bf I})=R_{G}(\Hat {\mbSigma}^{*}_{I}, {\bf I})
=R_{G}(\Hat {\mbSigma}_{GC}, {\bf I}) < R_{G}(\Hat {\mbSigma}_{S}, {\bf I})$ approximatly. We omit the details.

\indent Example 3 points out the fact that based on the orthogonal decomposition, the analytical difficulty to
find out the marginal expectations will be accompanied with the orthogonal decomposition due to the Selberg type
integral which involves the Vandermonde determinant. Similar difficulty will occur when to obtain the density
functions of $\mbox{min}_{1 \leq i \leq p}\{{l_i}\}$ and $\mbox{max}_{1 \leq i \leq p}\{{l_i}\}$ (for details
see Edelman, 1989). When $p, n \to \infty$ such that $\lim_{n\to \infty} p/n=y \in [0, 1]$, for a Wishart matrix
Edeiman (1989) proved  that the geometric means of
${\mathcal E}[\mbox{log}(\mbox{min}_{1 \leq i \leq p}\{{l_i}\}/n)]=\mbox{log}(1-\sqrt {y})^{2}+o(1)$
and ${\mathcal E}[\mbox{log}(\mbox{max}_{1 \leq i \leq p}\{{l_i}\}/n)]=\mbox{log}(1+\sqrt {y})^{2}+o(1)$,
respectively.

\indent With respect to Stein loss function, Stein (1956) made statistical inference focused on the special case
$\mbSigma={\bf I}$, it is sufficient enough due to the invariance consideration. For the main purpose of focusing
on statistical inference in this note, we adopt the same structure as Stein did. However, this will scarifice the
developement of distribution theory of arbitrary covariance $\mbSigma$. For this purpose, the elegant zonal
polynomials have been incorporated, we may refer the book of Muirhead (1982) for the details.

\indent The invariance nature of geodesic distance loss function suggests that it should be the way to deal with
state-of-the-art covariance matrix estimators for the interesting and timely large dimensional case. For the large
dimensional case, the sample size is required to be the same order of dimension. The emperical density of
eigenvalues of Wishart matrix ${\bf A}$ converges to the Marchenko-Paster law in the limit when $p, n \to \infty$
such that $c=p/n$ is fixed $0 \leq c  \leq 1.$ The geometric mean of this distribution is
$-1-[(1-y)\mbox{log}(1-y)]/y$, where $y=\lim_{n\to \infty} p/n$.

\vspace{0.3cm}
{\bf References}

\begin{enumerate}
{\small
\item \label{An}
{Anderson, T. W}. (2003). {\it An Introduction to Multivariate Statistical Analysis}, 3rd ed.. Wiley,
New York.

\item \label{As}
{Askey, R.} (1980). Some basic hypergeometric extensions of integrals of Selberg and Andrews. {\it SIAM,
J. Math. Anal.}, {\bf 11}, 938-951.

\item \label{Ea}
{Edelman, A}. (1989). Eigenvalues and condition numbers of random matrices. {\it Ph.D thesis},
Massachusetts institute of technology.

\item \label{JS}
{James, W. and Stein, C.} (1961). Estimation with quadratic loss. {\it Proc. Fourth Berleley Symp. Math.
Statist. Probab.}, {\bf 1}, 361-379. California Press, Berkeley, CA.

\item \label{Ki}
{Kiefer, J}. (1957). Invariance, minimax sequential estimation, and continuous time processes.
{\it Ann. Math. Statist.}, {\bf 28}, 573-601.

\item \label{An}
{Muirhead, R. J}. (1982). {\it Aspects of Multivariate Statistical Theory}. Wiley, New York.

\item \label{Te}
{Terras, A.} (1988). {\it Harmonic Analysis on Symmetric Spaces and Applications. II}. Springer, Berlin.





}
\end{enumerate}

\hfill


\end{document}